\documentclass[12pt,a4,amsmath]{article}
\usepackage[headings]{fullpage}
\usepackage{amssymb,epic,eepic,epsfig,amsbsy,amsmath}
\numberwithin{equation}{section} 
\usepackage{amscd}
\usepackage{times}
\def\HP{\mathop {\mathrm {HP}}\nolimits}
\def\K{\mathop {\mathrm {K}}\nolimits}
\def\E{\mathop {\mathrm {E}}\nolimits}
\def\I{\mathop {\mathbf {I}}\nolimits}
\def\Cyl{\mathop {\mathrm {Cyl}}\nolimits}

\def\Cone{\mathop{\mathrm {Cone}}\nolimits}

\def\K{\mathop {\mathrm {K}}\nolimits}

\def\Cyl{\mathop{\mathrm {Cyl}}\nolimits}
\def\ev{\mathop{\mathrm {ev}}\nolimits}
\def\iIm{\mathop{\mathrm {Im}}\nolimits}
\def\d{\mathop{\mathrm {d}}\nolimits}
\newtheorem{thm}{Theorem}[section]
\newtheorem{defn}[thm]{Definition}

\newtheorem{lem}[thm]{Lemma}
\newtheorem{cor}[thm]{Corollary}
\newtheorem{rem}[thm]{Remark}

\em
\begin{document}
\title{Category of Noncommutative CW complexes. III\footnote{The
work was supported in part by the National Project of Research in
Fundamental Sciences of Vietnam.}}
\author{Do Ngoc Diep} \maketitle
\begin{abstract}
We prove in this paper a noncommutative version of Leray Theorem
and then Leray-Serre Spectral Theorem for noncommutative Serre
fibrations: for NC Serre fibration there are converging spectral
sequences with $\E^2$ terms as $\E^2_{p,q} = \HP_p(A; \HP_q(B,A))
\Longrightarrow \HP_{p+q}(B)$ and $\E^2_{p,q} = \HP_p(A;\K_q(B,A))
\Longrightarrow \K_{p+q}(B)$.
\end{abstract}
\section*{Introduction}
The ideas of using spectral sequences to operator algebras was
started in \cite{diep1}. In that work the author constructed for
an arbitrary algebras some CCR- or CT- composition series and
tried to define the structure of the algebras through Busby
invariants of extensions. This reduced to the ideas of using
spectral sequences. However, the technique of spectral sequence
meets some difficulties, for example convergence problem. In this
paper we consider the NCCW structure in place of the composition
series in order to manage the spectral sequences in computing
$\HP$ or $\K$ groups for some NCCW complexes.

It is well-known that for every short exact sequence of ideal and
algebras \begin{equation}\CD 0 @>>> A @>>> B @>>> C @>>>
0\endCD\end{equation} where $A$ is some ideal in $B$ and $C \cong
B/A$ is the quotient algebra, there is a natural spectral sequence
converging to the cyclic periodic homology with $E_2$-term
\begin{equation}\E^2_{p,q} =  \HP_p(A;\HP_q(B,A)) \Longrightarrow
\HP_n(B),\mbox{ with } n=p+q\end{equation}
and a spectral sequence converging to the K-theory with $E_2$ term
\begin{equation}\E^2_{p,q} =  \HP_p(A;\K_q(B,A)) \Longrightarrow
\K_n(B),\mbox{ with } n=p+q.\end{equation}

It is natural to ask whether the condition of $A$ being an ideal
is necessary for this exact sequence. We observe \cite{diep3} that
this condition is indeed not necessarily to be satisfied. From
this we arrive to a more general condition of a noncommutative
Serre fibration as some homomorphism of algebras with the so
called HLP ({\it Homotopy Lifting Problem}). But we do restrict to
consider only the so called noncommutative CW-complexes. We proved
that up to homotopy we can change any homomorphism between
noncommutative CW-complexes and the CW-complexes them-selves to
have a noncommutative Serre fibration. We can then deduce the
spectral sequences for the cyclic periodic homology and K-theory
of an arbitrary NC Serre fibration and of an arbitrary morphism of
NC CW-complexes.
 It seems to the author that recently the notion of NC CW-complex
 is introduced , see e.g. \cite{cuntz} and here we  work with spectral
 sequences for general maps between NC CW-complexes.

Let us describe the contents of the paper. In Section 1 we
construct for noncommutative (NC) CW-complexes the Leray spectral
sequences. Next, in the section 2, we prove that for NC Serre
fibration there are spectral sequences with $\E^2$ terms as
\begin{equation}\E^2_{p,q} = \HP_p(A; \HP_q(B,A)) \Longrightarrow \HP_{p+q}(B)\end{equation}
and \begin{equation}\E^2_{p,q} = \HP_p(A;\K_q(B,A))
\Longrightarrow \K_{p+q}(B).\end{equation}

\section{Leray Spectral Sequence Theorem for NCCW}
We start with the notion of NCCW, introduced by S. Eilers, T. A.
Loring and G. K. Pedersen \cite{elp} and G. Pedersen
\cite{pedersen}. For reader's convenience, we repeat some
definitions form \cite{diep3}.

\begin{defn}\label{defn1.2} A dimension 0 NCCW complex {\rm is defined, following
\cite{pedersen} as a finite sum of C* algebras of finite linear
dimension, i.e. a sum of finite dimensional matrix algebras,}
\begin{equation}A_0 = \bigoplus_{k} \mathbf M_{n(k)}.\end{equation}
In dimension n, an NCCW complex  {\rm is defined
 as a sequence $\{A_0,A_1,\dots,A_n\}$ of C*-algebras $A_k$
 obtained each from the previous one by the pullback construction
 \begin{equation}\CD 0 @>>> \I^k_0F_k @>>> A_k @>\pi>> A_{k-1} @>>> 0\\
 @. @| @VV\rho_k V  @VV\sigma_k V @. \\
 0 @>>> \I_0^kF_k @>>> \I^kF_k @>\partial>> \mathbf S^{k-1}F_k @>>>
 0,\endCD \end{equation}
 where $F_k$ is some C*-algebra of finite
 linear dimension, $\partial$ the restriction morphism, $\sigma_k$
 the connecting morphism, $\rho_k$ the projection on the first
 coordinaes and $\pi$ the projection on the second coordinates in
 the presentation
 \begin{equation}
 A_k = \mathbf I^kF_k \bigoplus_{\mathbf S^{k-1}F_k}A_{k-1}\end{equation}
 }\end{defn}

\begin{defn}{\rm
We say that the morphism $f: A \longrightarrow B$ admits the so
called {\it Homotopy Lifting Property (HLP)} if for every algebra
$C$ and every morphism $\varphi : A \longrightarrow C$ such that
there is some morphism $\tilde{\varphi} : B \longrightarrow C$
satisfying $\varphi = \tilde{\varphi}\circ f$, and for every
homotopy $\varphi_t : A \longrightarrow C$, $\varphi_0 = \varphi$,
there exists a homotopy $\tilde{\varphi}_t : B \longrightarrow C$,
$\tilde\varphi_0 = \tilde\varphi$,  such that for every $t$,
$\varphi_t = \tilde{\varphi}_t \circ f$, i.e. the following
diagram is commutative }
\begin{center}\begin{equation}
\begin{picture}(60,60)
\put(0,0){A} \put(50,0){C} \put(0,50){B}
\put(10,5){\vector(1,0){35}} \put(5,10){\vector(0,1){35}}
\put(15,45){\vector(1,-1){30}} \put(10,25){f}
\put(37,25){$\tilde{\varphi}_t,\tilde{\varphi}_0=\tilde{\varphi}$}
\put(25,10){$\varphi_t$} \put(5,-10){$\varphi_0 =
\tilde{\varphi}_0 \circ f$}
\end{picture}
\end{equation}\end{center}
\end{defn}
\begin{defn}{\rm
A morphism of C*-algebras $f : A \longrightarrow B$ with HLP axiom
is called a {\it noncommutative Serre fibration} (NCSF).}
\end{defn}

\begin{thm}\label{thm1}
In the category of NCCW complexes, for every morphism $f : A
\longrightarrow B$, there is some homotopies $A \sim A'$, $B \sim
B'$ and a morphism $f' : A' \longrightarrow B'$ which is a NC
Serre fibration.
\end{thm}

Before prove this theorem we do recall \cite{diep1} the following
notions of noncommutative cylinder and noncommutative mapping
cone.

\begin{defn}[NC cone]{\rm For C*-algebras the} NC cone {\rm of $A$
is defined as the tensor product with } $\mathbf C_0((0,1])$, i.e.
\begin{equation}\Cone(A) := \mathbf C_0((0,1]) \otimes
A.\end{equation}
\end{defn}

\begin{defn}[NC suspension]{\rm For C*-algebras the} NC suspension
{\rm of $A$ is defined as the tensor product with } $\mathbf
C_0((0,1))$, i.e. \begin{equation}\mathbf S(A) := \mathbf
C_0((0,1))\otimes A.
\end{equation}
\end{defn}

\begin{rem}
If $A$ admits a NCCW complex structure, the same have the cone
$\Cone(A)$ of $A$ and the suspension $\mathbf S(A)$ of $A$.
\end{rem}

\begin{defn}[NC mapping cylinder]{\rm
Consider a map $f : A \to B$ between C*-algebras. The NC {\it
mapping cylinder} $\Cyl(f:A\to B)$ is defined by the pullback
diagram \cite{diep1}
\begin{equation}\CD \Cyl(f) @>pr_1 >> \mathbf C[0,1] \otimes A \\
@V{pr_2}VV  @VV{f\circ \ev(1)}V\\
B @>id>> B\endCD,\end{equation} where $\ev(1)$ is the map of
evaluation at the point $1\in [0,1]$. It can be also defined
directly as follows.
 In the algebra $\mathbf C(\mathbf I) \otimes A
\oplus B$ consider the closed two-sided ideal $\langle \{1\}
\otimes a - f(a), \forall a\in A\rangle $, generated by elements
of type $\{1\} \otimes a - f(a), \forall a\in A$. The quotient
algebra
\begin{equation}\Cyl(f)=\Cyl(f:A\to B) := \left(\mathbf C(\mathbf
I) \otimes A \oplus B\right)/ \langle \{1\} \otimes a - f(a),
\forall a\in A\rangle \end{equation} is called the} NC mapping
cylinder {\rm and denote it by } $\Cyl(f:A\to B)$.
\end{defn}

\begin{rem}
It is easy to show that $A$ is included in $\Cyl(f: A \to B)$ as
$\mathbf C\{0\} \otimes A \subset \Cyl(f:A\to B)$ and $B$ is
included in also $B \subset \Cyl(f:A\to B)$.
\end{rem}
\begin{defn}[NC mapping cone]{\rm The NC mapping cone
$\Cone(\varphi)$ is defined from the pullback diagram
\begin{equation}  \CD \Cone(\varphi) @>pr_1 >> \mathbf C_0(0,1] \otimes A \\
@V{pr_2}VV  @VV{f\circ \ev(1)}V\\
B @>id>> B\endCD,  \end{equation} where $\ev(1)$ is the map of
evaluation at the point $1\in [0,1]$. It can be also directly
defined as follows. In the algebra $\mathbf C((0,1]) \otimes A
\oplus B$ consider the closed two-sided ideal $\langle \{1\}
\otimes a - f(a), \forall a\in A\rangle $, generated by elements
of type $\{1\} \otimes a - f(a), \forall a\in A$. We define the}
mapping cone {\rm as the quotient algebra }
\begin{equation}\Cone(f)=\Cone(f:A\to B) := \left(\mathbf
C_0((0,1]) \otimes A \oplus B\right)/ \langle \{1\} \otimes a -
f(a), \forall a\in A\rangle .\end{equation}
\end{defn}
\begin{rem}
It is easy to show that $B$ is included in $\Cone(f: A \to B)$.
\end{rem}

\begin{thm}[Noncommutative Leray Theorem for HP-groups]
Let \begin{equation} \{0\}=A_{-1} \subseteq A_0 \subseteq \dots
\subseteq A_k = A\end{equation} be a NCCW complex structure of
$A$. Then for every $r \geq 0$, and all $p$, $q$, there are groups
$\E^r_{p,q}$ such that $\E^r_{p,q} = 0$ for all $p < 0$ or $p >
k$, and homomorphisms \begin{equation}\d^r_{p,q}: \E^r_{p,q} \to
\E^r_{p-r, q+r-1}\end{equation} such that $\d^r_{p-r,q+r-q}\circ
\d^r_{p,q} = 0$ and that
\begin{enumerate}
\item \begin{equation}\E^{r+1}_{p,q} = \ker \d^r_{p,q}/\iIm
\d^r_{p+r, q-r+1}.\end{equation} \item \begin{equation}\E^0_{p,q}=
C_q(A_p,A_{p-1}).\end{equation} \item $\d^0_{p,q}$ is coincided
with the boundary operator \begin{equation}\partial:
C_{p+q}(A_p,A_{p-1}) \to C_{p+q-1}(A_p,A_{p-1})\end{equation}
\item \begin{equation}\E^1_{p,q} =
\HP_{p+q}(A_p,A_{p-1}).\end{equation} \item $\d^1_{p,q}$ is
coincided with the homomorphism \begin{equation}\partial_{*}:
\HP_{p+q}(A_p,A_{p-1}) \to
\HP_{p+q-1}(A_{p-1},A_{p-2}),\end{equation} associated with the
triple $(A_p,A_{p-1},A_{p-2})$. \item
\begin{equation}\E^\infty_{p,q} = \frac{\iIm(\HP_{p+q}(A_p)\to
\HP_{p+q}(A))}{\iIm(\HP_{p+q}(A_{p-1})\to \HP_{p+q}(A))} =
\frac{{}_{(p)}\HP_{p+q}(A)}{{}_{(p-1)}\HP_{p+q}(A)}.\end{equation}
\end{enumerate}
\end{thm}
{\sc Proof}. The proof is quite similar to the same one in the
classical algebraic topology. There is no need of special
modification. \hfill$\Box$

\begin{thm}[Noncommutative Leray Theorem for K-groups]
Let \begin{equation}\{0\}= A_{-1} \subseteq A_0 \subseteq \dots
\subseteq A_k = A\end{equation} be NCCW complex structure of $A$.
Then for every $r \geq 0$, and all $p$, $q$, there are groups
$\E^r_{p,q}$ such that $\E^r_{p,q} = 0$ for all $p < 0$ or $p >
k$, and homomorphisms \begin{equation}\d^r_{p,q}: E^r_{p,q} \to
\E^r_{p-r, q+r-1}\end{equation} such that $\d^r_{p-r,q+r-q}\circ
\d^r_{p,q} = 0$ and that
\begin{enumerate}
\item \begin{equation}\E^{r+1}_{p,q} = \ker \d^r_{p,q}/\iIm
\d^r_{p+r, q-r+1}.\end{equation} \item \begin{equation}\E^0_{p,q}=
C(A_p,A_{p-1}).\end{equation} \item $\d^0_{p,q}$ is coincided with
the boundary operator \begin{equation}\partial:
C_{p+q}(A_p,A_{p-1}) \to C_{p+q-1}(A_p,A_{p-1}).\end{equation}
\item \begin{equation}\E^1_{p,q} =
\K_{p+q}(A_p,A_{p-1})\end{equation} \item $\d^1_{p,q}$ is
coincided with the homomorphism \begin{equation}\partial_{*}:
\K_{p+q}(A_p,A_{p-1}) \to
\K_{p+q-1}(A_{p-1},A_{p-2}),\end{equation} associated with the
triple $(A_p,A_{p-1},A_{p-2})$. \item
\begin{equation}\E^\infty_{p,q} = \frac{\iIm(\K_{p+q}(A_p)\to
\K_{p+q}(A))}{\iIm(\K_{p+q}(A_{p-1})\to \K_{p+q}(A))} =
\frac{{}_{(p)}\K_{p+q}(A)}{{}_{(p-1)}\K_{p+q}(A)}.\end{equation}
\end{enumerate}
\end{thm}
{\sc Proof}. The proof is quite similar to the same one in the
classical algebraic topology. There is no need of special
modification. \hfill$\Box$

\section{Spectral Sequences for Noncommutative Serre Fibrations}

\begin{defn}
Let $A\hookrightarrow B$ be a Serre fibration.  The system of all
homotopy equivalences \begin{equation}i: (B:A) := \Pi(A,Cone(B,A))
\to B\end{equation} is called the {\bf system of  local
coefficients}. We say that a system of local coefficients is {\bf
simple} iff the K-theory and cyclic theory of $(B:A)$ is
independent of the choice of $i$.
\end{defn}

From now on we consider only  NC Serre fibrations with simple
system of local coefficients.

\begin{lem} For the Serre fibration $f: A \hookrightarrow B$ such that the
spectrum $\hat{B} \setminus \hat{A}$ has a cellular decomposition
\begin{equation}\hat{B} \setminus \hat{A}=\bigcup_\alpha (\mathbf
D^p_\alpha,\mathbf S^{p-1}_\alpha),\quad \mathbf
D^p_\alpha=\mathbf D^p, \mathbf S^{p-1}_\alpha=\mathbf
S^{p-1}\end{equation} i.e. \begin{equation}B = A
\bigoplus_{\mathbf S^{p-1}F_p} \mathbf I^p F_p\end{equation} for
some finite set $F_p$, we have
\begin{equation}A \otimes \mathbf C([0,1]) \oplus B /(A \times
\{0\} \oplus f(A) \oplus A \otimes \mathbf C(\{1\})) \cong A \otimes
\mathcal K \otimes \mathbf C(\hat{B}\setminus \hat{A})\end{equation}
\end{lem}

\begin{thm}
For an arbitrary NC Serre fibration $A \hookrightarrow B$   with
simple system of local coefficients, we have
\begin{enumerate}
\item \begin{equation}\E^1_{p,q} = \mathcal
C_p(A;\HP_q(B,A)),\end{equation} \item \begin{equation}\d^1_{p,q}
=
\partial : \mathcal C_p(A;\HP_q(B,A)) \to \mathcal
C_{p-1}(A;\HP_q(B,A)),\end{equation} \item
\begin{equation}\E^2_{p,q} = \HP_p(A;\HP_q(B,A)).\end{equation}
\end{enumerate}
\end{thm}
{\sc Proof}. In the category of NCCW complexes, we have the fact
that the graded de Rham cohomology of a cell is isomorphic to the
periodic cyclic homology of the corresponding NC cell. And
therefore
$$\E^1_{p,q} = \frac{\ker\{\HP_q(A_p,A_{p-1}\to \HP_q(A_{p-1},A_{p-2}))
\}}{\ker\{\HP_q(A_p,A_{p-1}\to \HP_q(A_{p},A_{p-1})) \}} =$$
$$= \HP_{p+q}(A_p,A_{p-1}) \cong \HP(\bigvee_{\mbox{p-dimensional cells}}
C(\mathbf S^p)\otimes \mathcal K \otimes A) \cong $$ $$\cong
\bigvee_{\mbox{p-dimensional cells}}\HP_{p+q}(\mathbf  C(\mathbf
D^p)\otimes \mathcal K \otimes B,C(\mathbf S^{p-1}) \otimes
\mathcal K \otimes A)  \cong$$ \begin{equation}\cong
\bigvee_{\mbox{p-dimensional cells}}
 \HP_{q}(\mathbf C)\otimes \HP_q(A) \cong \mathcal C_q(B,A; \HP_q(\mathbf C)).\end{equation}
It is easy to see that
\begin{equation}\tilde{\d}^1_{p,q} = \delta : \mathcal C_p(A,\tilde{\HP}_q(B,A) \to
\mathcal C_{p+1}(A;\tilde{\HP}_q(B,A))\end{equation} and therefore
\begin{equation}\tilde{\E}^2_{p,q}A= \HP(A;\tilde{\HP}_q(B,A)).\end{equation}
\hfill$\square$

\begin{cor}
In the category of NCCW complexes for an arbitrary map $A \to B$,
we have
\begin{enumerate}
\item \begin{equation}\E^1_{p,q} = \mathcal
C_p(A,\HP_q(B,A)),\end{equation} \item \begin{equation}\d^1_{p,q}
=
\partial : \mathcal C_p(A;\HP_q(B,A)) \to \mathcal
C_{p-1}(A,\HP_q(B,A)),\end{equation} \item
\begin{equation}\E^2_{p,q} = \HP_p(A;\HP_q(B,A)).\end{equation}
\end{enumerate}
\end{cor}

By analogy we have also the same results for K-groups
\begin{thm}
For an arbitrary NC Serre fibration $A \hookrightarrow B$ with
simple system of local coefficients, we have
\begin{enumerate}
\item \begin{equation}\E^1_{p,q} = \mathcal
C_p(A,\K_q(B,A)),\end{equation} \item \begin{equation}\d^1_{p,q} =
\partial : \mathcal C_p(A,\K_q(B,A)) \to \mathcal
C_{p-1}(A,\K_q(B,A)),\end{equation} \item
\begin{equation}\E^2_{p,q} = \HP_p(A;\K_q(B,A)).\end{equation}
\end{enumerate}
\end{thm}
{\sc Proof}. In the category of NCCW complexes, we have the fact
that the graded de  Rham cohomology of a cell is isomorphic to the
periodic cyclic homology of the corresponding NC cell. And
therefore
$$\E^1_{p,q} = \frac{\ker\{\K(A_p,A_{p-1}\to \K(A_{p-1},A_{p-2})) \}}{\ker\{\K(A_p,A_{p-1}\to \K(A_{p},A_{p-1})) \}} =$$
$$= \K_{p+q}(A_p,A_{p-1}) \cong \K(\bigvee_{\mbox{p-dimensional cells}}
C(\mathbf D^p)\otimes \mathcal K \otimes A) \cong $$ $$\cong
\bigvee_{\mbox{p-dimensional cells}}\K_{p+q}(\mathbf C(\mathbf
D^p)\otimes \mathcal K \otimes A,\mathbf S^{p-1}\otimes \mathcal K
\otimes A) \cong$$ \begin{equation}\cong
\bigvee_{\mbox{p-dimensional cells}} \K_{q}(A) \cong \mathcal
C_q(A; \K_q(B,A)).\end{equation} It is easy to see that
\begin{equation}\tilde{\d}^1_{p,q} = \delta : \mathcal C_p(A;\tilde{\K}_q(B,A)) \to \mathcal C_{p+1}(A;\tilde{\K}_q(B,A))\end{equation} and therefore
\begin{equation}\tilde{\E}^2_{p,q} = \HP_p(A;\tilde{\K}_q(B,A)).\end{equation}
\hfill$\square$

\begin{cor}
In the category of NCCW complexes for an arbitrary morphism $A \to
B$, we have
\begin{enumerate}
\item \begin{equation}\E^1_{p,q} = \mathcal
C_p(A;\K_q(B,A)),\end{equation} \item \begin{equation} \d^1_{p,q}
=
\partial : \mathcal C_p(A;\K_q(B,A)) \to \mathcal
C_{p-1}(A;\K_q(B,A)),\end{equation} \item
\begin{equation}\E^2_{p,q} = \HP_p(A;\K_q(B,A)).\end{equation}
\end{enumerate}
\end{cor}
\vskip 0.5cm The following remark was suggested by C. Schochet
\begin{rem}{\rm \label{rem2.6}
These theorem are true because both K-theory and cyclic theory
have very strong excision properties (analogous in spaces to the
property that given a space $X$ and a closed subspace $Y$ with $f:
Y \to X$ the inclusion, that the natural map $\Cone(f) \to X/Y$
induces an isomorphism on K-theory and cyclic theory (but this
fails for many other theories.)) }\end{rem}

\begin{rem}{\rm
In a subsequent paper, we shall  prove that for the induced
noncommutative Serre fibration with simple system of local
coefficients there is also noncommutative analog of the
Eilenberg-Moore spectral sequences, which are collapsing at $\E^2$
terms if the system of local coefficients are simple. }\end{rem}

\section*{Acknowledgments}
The author would like to thank Professor Claude Schochet for very
useful inspirational  e-Mail messages and in particular for some
suggestion concerning terminology, Remak \ref{rem2.6} and  for the
reference \cite{schochet}.

{\sc Institute of Mathematics, Vietnam Academy of Sciences and Technology, 18 Hoang Quoc Viet Road, Cau Giay district , 10307 Hanoi, Vietnam}\\
{\tt e-Mail: dndiep@math.ac.vn}
\end{document}